\documentclass[11pt, a4paper]{article}
\usepackage{array}
\usepackage{multirow}
\makeatletter
\newcommand{\thickhline}{%
    \noalign {\ifnum 0=`}\fi \hrule height 1pt
    \futurelet \reserved@a \@xhline
}
\newcolumntype{"}{@{\hskip\tabcolsep\vrule width 1pt\hskip\tabcolsep}}
\makeatother

\usepackage{graphics}
\usepackage{graphicx}
\usepackage{epsfig}
\usepackage{subfigure}
\usepackage{amssymb}
\usepackage{amsthm}
\usepackage{mathrsfs}
\usepackage{hhline}
\usepackage{longtable}
\usepackage{amsmath}
\usepackage{array}
\usepackage{float}
\usepackage{multirow}
\usepackage{cite}
\usepackage{enumerate}
\usepackage{xcolor}

\def\ms{\medskip}
\def\nt{\noindent}

\newtheoremstyle{de}
  {10pt}          
  {10pt}  
  {\rm}  
  {}
  {\bf}  
  {. }    
  { }    
  {}     
\theoremstyle{de}

\newtheorem{example}{Example}[section]

\newtheoremstyle{theorem}
  {10pt}          
  {10pt}  
  {\it}  
  {}
  {\bf}  
  {. }    
  { }    
  {}     
\theoremstyle{theorem}
\newtheorem{theorem}{Theorem}[section]
\newtheorem{lemma}[theorem]{Lemma}

\numberwithin{equation}{section}

\graphicspath{%
    {converted_graphics/}
    {/}
}

\begingroup\makeatletter\ifx\SetFigFont\undefined%
\gdef\SetFigFont#1#2#3#4#5{%
  \reset@font\fontsize{#1}{#2pt}%
  \fontfamily{#3}\fontseries{#4}\fontshape{#5}%
  \selectfont}%
\fi\endgroup%

\definecolor{vividviolet}{rgb}{0.62, 0.0, 1.0}

\vspace{5cm}
\begin{document}

{
\title { {\Large \sc {\bf{On local antimagic chromatic number of a corona product graph}}}}
\author{{\sc Gee-Choon Lau{$^a$}, M. Nalliah$^{b,}$\footnote{Corresponding author.} }\\
	{\footnotesize $^a$Faculty of Computer \& Mathematical Sciences}\\
	{\footnotesize Universiti Teknologi MARA (Segamat Campus)}\\
	{\footnotesize 85000, Johor, Malaysia.}\\
	{\footnotesize $^b$School of Advanced Sciences}\\
	{\footnotesize Department of  Mathematics }\\
	{\footnotesize  Vellore Institute of Technology}\\ 
	{\footnotesize  Vellore-632 014,  India.}\\ 
	{\footnotesize e-mail: {\ttfamily geeclau@yahoo.com,  nalliahklu@gmail.com }}}
\date{}
}

\maketitle

\begin{abstract}
In this paper, we provide a correct proof for the lower bounds of the local antimagic chromatic number of the corona product of friendship and fan graphs with null graph respectively as in [On local antimagic vertex coloring of corona products related to friendship and fan graph, {\it Indon. J. Combin.}, 5(2) (2021) 110--121]. Consequently, we obtained a sharp lower bound that gives the exact local antimagic chromatic number of the corona product of friendship and null graph. 
\end{abstract}



\section{Introduction}

Let $G=(V,E)$ be a finite, undirected graph with neither loops nor multiple edges. The order and size of $G$ are denoted by $|V(G)|=p$ and $|E(G)|=q$ respectively. The {\it friendship graph} $f_n$ $(n\ge 2)$ is a graph which consists of $n$ triangles with a common vertex. The {\it fan graph} $F_n$ $(n\ge 2)$ is obtained by joining a new vertex to every vertex of a path $P_n$. The \textit{corona product} of two graphs $G$ and $H$ is the graph $G\circ H$ obtained by taking one copy of $G$ along with $|V(G)|$ copies of $H$, and join the $i$-th vertex of $G$ to every vertex of the $i$-th copy of $H$, where $1 \leq i \leq |V(G)|$. For integers $a < b$, let $[a,b]=\{a,a+1,\ldots,b\}$.
For graph-theoretic terminology, we refer to Chartrand and Lesniak \cite{cl}. 

\ms\nt Hartsfield and Ringel's \cite{HR} introduced the concept of antimagic labeling of a graph. For a graph $G$, let $f:E(G)\to \{1, 2, . . . , q\}$ be a bijection. For each vertex $u\in V(G)$, the weight $w(u)=\sum_{e\in E(u)}f(e)$, where $E(u)$ is the set of edges incident to $u$. If $w(u)\neq w(v)$ for any two distinct vertices $u$ and $v\in V(G)$, then $f$ is called an antimagic labeling of $G$.  Hartsfield and Ringel conjectured that every connected graph with at least three vertices admits antimagic labeling \cite{HR}. Interested readers can refer to \cite{TEC, JAG}.
 
\ms\nt  Arumugam et al. in \cite{SA}, and independently, Bensmail et al. in \cite{JB}, posed a new  definition as a relaxation of the notion of antimagic labeling. They called a bijection $f:E\rightarrow \left\{1, 2, . . . , \left| E\right|\right\}$ a {\it local antimagic labeling} of $G$ if for any two adjacent vertices $u$ and $v$ in $V (G),$ the condition  $w(u) \neq w(v)$ holds. Based on this notion, Arumugam et al. then introduced a new graph  coloring parameter. 
%
%
Let $f$ be a local antimagic labeling of $G$ be a connected graph $G$. The assignment of $w(u)$ to $u$ for each vertex $u\in V(G)$ induces naturally a proper vertex coloring of $G$ which is called a {\it local antimagic vertex coloring} of $G.$ The {\it local antimagic chromatic number}, denoted $\chi_{la}(G)$, is the minimum number of colors taken over all local antimagic colorings of $G$~\cite{SA}.
 
 


\ms\nt Arumugam et al. \cite{SA2} obtained the local antimagic chromatic number for the graph $G \circ O_m$, where $G$ is a path, cycle or complete graph and $O_m$ is the null graph of order $m\ge 1$.

\begin{theorem}\rm{\cite{SA2}\label{sa2}}
	Let $m\geq 2$, then $\chi_{la} \left(C_3 \circ O_m\right) =3m+3$ except $\chi_{la}(C_3\circ O_1)=5$.    
\end{theorem} 

\begin{theorem}\rm{\cite{SA2}\label{sa2a}}
For $n\ge 2$, $\chi_{la} \left(K_n \circ K_1\right)=2n-1$.
\end{theorem} 

\nt In~\cite{Himami+S}, the authors studied $\chi_{la}(f_n\circ O_m)$ and $\chi_{la}(F_n \circ O_m)$ for $n\ge 2$ and $m\ge1$. We note that there are inconsistencies in the notations of $f_n$ and $F_n$ used. They proved that $\chi_{la}(f_n\circ O_m) \le m(2n+1)+3$ and $\chi_{la}(F_n) \le m(n+1)+3$ by providing a correct local antimagic labeling respectively. However, there are gaps in proving that $\chi_{la}(f_n\circ O_m) \ge m(2n+1)+3$ and $\chi_{la}(F_n) \ge m(n+1)+3$. Motivated by this, we shall first provide correct arguments to the proofs of the lower bounds. Consequently, we showed that $\chi_{la}(f_n \circ O_m) = m(2n+1)+2$ for $n\ge 2, m=1$. Interested readers may refer to~\cite{LHS, LSN, LSN-pendant, LSS} for local antimagic chromatic number of graphs with pendant edges. 

\section{Lower bounds of $\chi_{la}(f_n\circ O_m)$ and  $\chi_{la}(F_n\circ O_m)$}


\begin{lemma}\label{lem-fnOm} For $n\ge 2, m\ge 1$, $\chi_{la}(f_n\circ O_m) \ge m(2n+1)+3$ except $\chi_{la}(f_n\circ O_1) \ge m(2n+1)+2$. \end{lemma}

\begin{proof} Let $G=f_n\circ O_m$ with $V(G)=\{x,u_i,v_i,x_j,u^i_j,v^i_j\,|\,1\le i\le n, 1\le j\le m\}$ and $E(G)=\{xx_j, xu_i, xv_i, u_iv_i, u_iu^i_j, v_iv^i_j\,|\,1\le i\le n, 1\le j\le m\}$. Clearly, $|E(G)|=q=m(2n+1)+3n$.

\ms\nt Suppose $f: E(G)\to [1,q]$ is a local antimagic labeling of $G$. Clearly, all the $m(2n+1)$ pendant vertices must have distinct induced vertex colors that are at most $q$. Morever, $w(x)\ge 1+2+\cdots+(2n+m) = (2n+m)(2n+m+1)/2=s$. Now, $2s-2q = (2n+m+1)^2 + (2n+m+1) - 6n - 2m(2n+1) = 4n^2+m^2+m+1 > 0$. Thus, $w(x) > q$. Therefore, $\chi_{la}(G)\ge m(2n+1)+1$. Without loss of generality, we consider the following 3 cases. 

\ms\nt{\bf Case 1.} $f(u_1v_1)=q$. In this case, $w(u_1)\ne w(v_1)\ne w(x) > q$ so that $\chi_{la}(G)\ge m(2n+1)+3$.

\ms\nt{\bf Case 2.} $f(xu_1)=q$ or $f(u_1u^1_1)=q$. In this case, $w(u_1) \ne w(x) > q$ so that $\chi_{la}(G)\ge m(2n+1)+2$. Suppose equality holds. Clearly, for each $i\in[1,n]$, at most one of $u_i, v_i$ has induced vertex color $q$. So, there are at most $n$ vertices in $\{u_i,v_i\}$ with induced vertex color $q$. The sum of these $n$ induced vertex colors is at least $1+2+\cdots + n(m+2) = n(m+2)[n(m+2)+1]$ and at most $nq = n[3n+m(n+1)]$. Since $n\ge 2$, it is easy to check that $n(m+2)[n(m+2)+1] - n[3n+m(n+1)] = 2n^2(m+1)+n+\frac{1}{2}mn(mn+1) - [3n^2+mn(2n+1)] > 0$ if and only if $m > 1$. Consequently, $\chi_{la}(G)\ge m(n+1)+2$ if $m=1$, and $\chi_{la}(G)\ge m(n+1)+3$ if $m\ge 2$.

\ms\nt{\bf Case 3.} $f(xx_1)=q$. In this case, $w(x_1)=q$ and $w(u^i_j), w(v^i_j), w(x_j) < q$ $(x_j\ne x_1)$ so that $\chi_{la}(G)\ge m(2n+1)+1$. Suppose $w(v_i) < w(u_i) \le q$ for $1\le i\le n$, then $\sum^n_{i=1}[w(u_i)+w(v_i)]$ is at most $n(2q-1)$ and at least $1+2+\cdots + n(2m+3)=n(2m+3)[n(2m+3)+1]/2$. Now,

\begin{align*}
&n(2m+3)[n(2m+3)+1] - 2n(2q-1)\\
&\quad=n(2m+3)[n(2m+3)+1] - 2n[2m(2n+1)+6n-1]\\
&\quad=4m^2n^2 + 4mn^2 - 2mn - 3n^2 + 5n > 0.
\end{align*}

\nt Thus, we may assume $w(u_1) > q$. Since $w(u_1)\ne w(x)$, we have $\chi_{la}(G)\ge m(2n+1)+2$. Suppose equality holds. By an argument similar to that in Case 2, we have $\chi_{la}(G)\ge m(2n+1)+ 2$ if $m=1$ and $\chi_{la}(G)\ge m(2n+1)+3$ if $m\ge 2$. 
\end{proof}

\nt Note that $F_2 \circ O_m = C_3 \circ O_m$, we next consider $F_n \circ O_m, n\ge 3, m\ge 1$.

\begin{lemma}\label{lem-FnOm} For $n\ge 3, m\ge 1$, $\chi_{la}(F_n \circ O_m) \ge m(n+1) + 3$. \end{lemma} 

\begin{proof} Let $G = F_n \circ O_m$ with $V(G) = \{x, x_j, v_i, v^i_j\,|\, 1\le i\le n, 1\le j\le m\}$ and $E(G) = \{xx_j, xv_i,  v_iv^i_j\,|\, 1\le i\le n, 1\le j\le m\} \cup \{v_iv_{i+1}\,|\, 1\le i\le n-1\}$. Clearly, $|E(G)| = m(n+1) + 2n-1 = q$. 


\ms\nt Let $f$ be a local antimagic labeling of $G$ that induces $\chi_{la}(G)$ distinct vertex colors. Clearly, all the $m(n+1)$ pendant vertices must have distinct induced vertex colors that are at most $q$. Moreover, $w(x) \ge 1+2+\cdots + (m+n)(m+n+1)/2 = s$. Now $2s - 2q = (m+n)(m+n+1) - 2[m(n+1)+2n-1] = m^2-m + n^2-3n+1 > 0$ for $n\ge 3$. Thus, $w(x) > q$ and $\chi_{la}(G)\ge m(n+1)+1$. Without loss of generality, we consider the following  cases.

\ms\nt{\bf Case 1.} $f(v_1v_2)=q$ or $f(v_2v_3)=q$ if $n\ge 4$. In this case, $w(x) \ne(v_1)\ne w(v_2) > q$. Thus, $\chi_{la}(G)\ge m(n+1)+3$.

\ms\nt{\bf Case 2.} $f(xv_1)=q$ (or $f(xv_2)=q$). In this case, $w(x)\ne w(v_1) > q$ (or $w(x)\ne w(v_2) > q$). Thus, $\chi_{la}(G)\ge m(n+1)+2$. Suppose equality holds. Note that if $w(v_i) > q$ for $3\le i\le n$, then $w(v_i) = w(v_1)$. Moreover, $w(v_i)\ne w(v_{i+1})$ for $1\le i\le n-1$. Suppose there are $r\ge 1$ vertices  in $\{v_i\,|\,1\le 1\le n\}$ with induced vertex color larger than $q$, then there are $n-r\ge 1$ vertices in $\{v_i\,|\,1\le 1\le n\}$ with induced vertex color at most $q$. These $n-r$ vertices are incident to a total of $(m+2)n-1 - r(m+1) = (m+1)(n-r)+n-1$. Therefore, their edge labels sum under $f$ is at most $(n-r)q$. However, the sum is at least $S=1+2+\cdots+[(m+1)(n-r)+n-1]=\frac{1}{2}[(m+1)(n-r)+n-1][(m+1)(n-r)+n]$. Note that $n-r\ge n/2$. Thus, $-r\ge -n/2$ and
\begin{align*}
&2S-2(n-r)q\\
&= [(m+1)(n-r)+n]^2-[(m+1)(n-r)+n]-2(n-r)[m(n+1)+2n]\\
&= (m+1)^2(n-r)^2+(2n-1)(m+1)(n-r)+n^2-n-2(n-r)(mn+m+2n)\\
&= (n-r)[m^2(n-r)+2m(n-r)-3m-n-r-1]+n^2-n\\
&\ge (n-r)[m^2(n-r)+2m(n-r)-3m-\frac{3n}{2}-1]+n^2-n\\
&\ge \frac{n}{2}\bigg[(m^2+2m)(\frac{n}{2})-3m-\frac{3n}{2}-1\bigg]+n^2-n\\
&\ge \frac{n}{2}\bigg[\frac{m^2n}{2}-\frac{3n}{2}-1+2n-2\bigg]\\
&= \frac{n}{2}\bigg[\frac{m^2n}{2}+\frac{n}{2}-3\bigg] > 0 \mbox{ except $n=3,m=1$,}
\end{align*}
\nt contradicting  $S \le (n-r)q$ for all $(n,m)\ne (3,1)$.  Now, consider $G=F_3\circ O_1$ that has $q=9$. If $G$ admits a local antimagic labeling that induces $6$ distinct vertex colors, then $w(v_1)=w(v_3)\le 9$. Since $v_1$ and $v_3$ are incident to 6 different edges, their total label sum is at least 21 so that $w(v_1)=w(v_3)\ge 11$, a contradiction. Therefore, $\chi_{la}(G)\ge m(n+1)+3$. 

\ms\nt{\bf Case 3.} $f(v_1v^1_1)=q$ (or $f(v_2v^2_1)=q$). In this case, $w(v_1) \ne w(x) > q$ (or $w(v_2)\ne w(x) > q$). Thus, $\chi_{la}(G)\ge m(n+1)+2$. Suppose equality holds. By an argument similar to Case 2, we have the same contradiction.   \end{proof}



\section{$\chi_{la}(f_n\circ O_1$)}

In~\cite{Himami+S}, the authors obtained local antimagic labelings that correctly show that $\chi_{la}(f_n\circ O_m)\le m(2n+1)+3$ and $\chi_{la}(F_n\circ O_m)\le m(n+1)+3$. By Lemma~\ref{lem-fnOm}, we shall next show that $\chi_{la}(f_n\circ O_1)=2n+3$.

\begin{theorem}\label{thm-fnO1} For $n\ge 2$, $\chi_{la}(f_n\circ O_1)=2n+3$. \end{theorem} 

\begin{proof} Let $G=f_n\circ O_1$ with $V(G)$ and $E(G)$ as defined in the proof of Lemma~\ref{lem-fnOm}. Suffice to define a bijection $f: E(G)\to[1,5n+1]$ that induces $2n+3$ distinct induced vertex colors. We shall use labeling matrices to describe the labeling of all the edges of $f_n\circ O_1$.   


\ms\nt Suppose $n$ is odd. We first define $f(xx_1)=5n+1$. We now arrange integers in $[2n+1,5n]$ as a $3\times n$ matrix as follows:
\begin{enumerate}[(1).]
  \item In row 1, assign $4n+(i+1)/2$ to column $i$ if $i=1,3,5\ldots,n$; assign $(9n+1)/2+i/2$ if $i=2,4,6,\ldots,n-1$. We have used integers in $[4n+1,5n]$.
  \item In row 2, assign $(7n+1)/2+(i-1)/2$ to column $i$ if $i=1,3,5\ldots,n$; assign $3n+i/2$ if $i=2,4,6,\ldots,n-1$. We have used integers in $[3n+1,4n]$.
  \item In row 3, assign $3n+1-i$ to column $1\le i\le n$. We have used integers in $[2n+1,3n]$.
\end{enumerate} 

\nt  The resulting matrix is given by the following table:
\begin{center}
\nt Table 1: assignment of integers in $[2n+2,5n+1]$\\
$\setlength{\extrarowheight}{2pt}
\begin{array}{|c|c|c|c|c|c|c|c|c|}\hline
4n+1 & \frac{9n+3}{2} & 4n+2 & \frac{9n+5}{2} & \cdots & 5n-1 & \frac{9n-1}{2} & 5n & \frac{9n+1}{2} \\[2pt]\hline
\frac{7n+1}{2} & 3n+1 & \frac{7n+3}{2} & 3n+2 & \cdots &  \frac{7n-3}{2} & 4n-1 &  \frac{7n-1}{2}  & 4n \\[2pt]\hline
3n & 3n-1 & 3n-2 & 3n-3 & \cdots & 2n+4 & 2n+3 & 2n+2 & 2n+1\\\hline
\end{array}$
\end{center}

\ms\nt We next arrange integers in $[1,3n]$ as a $3\times n$ matrix as follows:
\begin{enumerate}[(1).]
  \item In row 1, assign $3n+1-i$ to column $1\le i\le n$. We have used integers in $[2n+1,3n]$.
  \item In row 2, assign $(3n+1)/2+(i-1)/2$ to column $i$ if $i=1,3,5,\ldots,n$; assign $n+i/2$ to column $i$ if $i=2,4,6,\ldots,n-1$. We have used integers in $[n+1,2n]$.
  \item In row 3, assign $(i+1)/2$ to column $i$ if $i=1,3,5,\ldots,n$; assign $(n+1)/2+i/2$ to column $i$ if $i=2,4,6,\ldots,n-1$. We have used integers in $[1,n]$.
\end{enumerate}

\nt  The resulting matrix is given by the following table:
\begin{center}
\nt Table 2: assignment of integers in $[2n+2,5n+1]$\\
$\setlength{\extrarowheight}{2pt}
\begin{array}{|c|c|c|c|c|c|c|c|c|}\hline
3n & 3n-1 & 3n-2 & 3n-3 & \cdots & 2n+4 & 2n+3 & 2n+2 & 2n+1\\\hline
\frac{3n+1}{2}& n+1 & \frac{3n+3}{2} & n+2 & \cdots & \frac{3n-3}{2} & 2n-1 & \frac{3n-1}{2} & 2n \\[2pt]\hline
1 & \frac{n+3}{2} & 2 & \frac{n+5}{2} & \cdots & n-1 & \frac{n-1}{2}   & n &  \frac{n+1}{2} \\[2pt]\hline
\end{array}$
\end{center}

\ms\nt For $1\le k\le 3$, $1\le i\le n$, let $a_{k,i}$ be the $(k,i)$-entry of Table 1, and $b_{k,i}$ be the $(k,i)$-entry of Table 2. Note that $b_{1,i} = a_{3,i}$. Define $f(u_iu^i_1) = a_{1,i}$, $f(xu_i) = a_{2,i}$, $f(u_iv_i) = a_{3,i}$, $f(xv_i) = b_{2,i}$ and $f(v_iv^i_1) = b_{3,i}$. It is obvious that $f$ is a bijective function.

\ms\nt Now, column sum of each column of Table 1 is $(21n+3)/2$. Thus, $w(u_i)=(21n+3)/2$ and $w(u^i_1)\in[4n+1,5n]$ for $1\le i\le n$. Similarly, the column sum of each column of Table 2 is $(9n+3)/2$. Thus, $w(v_i)=(9n+3)/2$ and $w(v^i_1)\in[1,n]$ for $1\le i\le n$. Moreover, $w(x) = (n+1) + \cdots + (2n) + (3n+1) + \cdots + 4n + (5n+1) = (n+1)(5n+1)$.  Clearly, $w(x) \ne w(u^i_1) \ne  w(u_i) \ne w(v^i_1) \ne w(x_1) = 5n+1$ for $1\le i\le n$. Note that $4n+1\le w(v_i) = (9n+3)/2 \le 5n+1$ odd for $n\ge 3$. Therefore, $f$ is a local antimagic labeling that induces $2n+3$ distinct vertex colors. Consequently, $\chi_{la}(f_n \circ O_1) = 2n+3$ for odd $n\ge 3$.

\ms\nt We now consider even $n\ge 2$. Figures~\ref{fig:n2m1} and~\ref{fig:n4m1} show that $\chi_{la}(f_2\circ O_1) = 7$ and $\chi_{la}(f_4\circ O_1) = 11$. 


\begin{figure}[h] 
  \centering
  \includegraphics[bb=3 3 266 131,width=3.1in,height=1.2in,keepaspectratio]{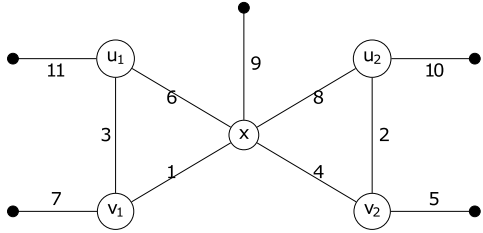}
  \caption{$\chi_{la}(f_2\circ O_1) = 7$ with induced vertex colors in $\{7,5,9,10,11,20,28\}$}
  \label{fig:n2m1}
\end{figure}

\begin{figure}[h] 
  \centering
  \includegraphics[bb=17 18 375 335,width=3.2in,height=3.2in,keepaspectratio]{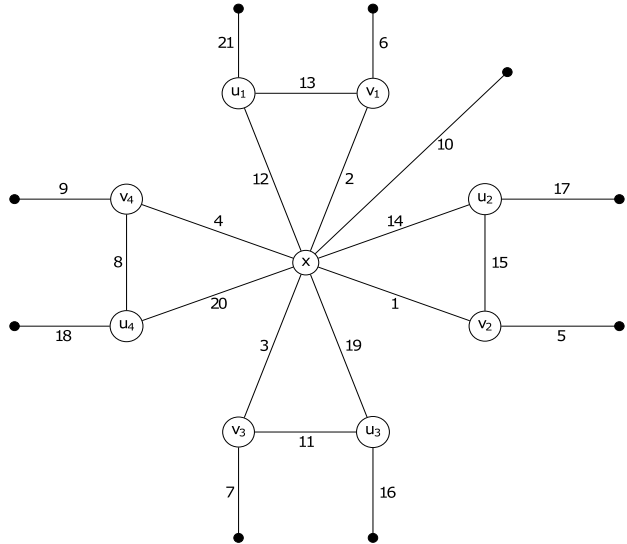}
  \caption{$\chi_{la}(f_4\circ O_1) = 11$ with induced vertex colors in $\{5,6,7,9,10,16,17,18,21,46,85\}$}
  \label{fig:n4m1}
\end{figure}

\ms\nt Consider $n\ge 6$. We first define $f(xx_1) = 3n+3$, $f(u_nv_n)=1$, $f(u_nu^n_1)=2n+2$, $f(xu_n)=2n$, $f(v_nv^n_1)=2n+3$ and $f(xv_n)=2n+1$. We now have $w(x_1)=3n+3$, $w(u_n)=4n+3$, $w(u^n_1)=2n+2$, $w(v_n)=4n+5$ and $w(v^n_1)=2n+3$. We now consider the remaining integers in $[2,2n-1]\cup[2n+4,3n+2]\cup[3n+4,5n+1]$. 

\ms\nt We now arrange integers in $[2n+4,3n+2]\cup[3n+4,5n+1]$ as a $3\times (n-1)$ matrix as follows:
\begin{enumerate}[(1).]
  \item In row 1, assign $4n+3+(i-1)/2$ to column $i$ if $i=1,3,5\ldots,n-1$; assign $9n/2+2+i/2$ if $i=2,4,6,\ldots,n-2$. We have used integers in $[4n+3,5n+1]$.
  \item In row 2, assign $7n/2+3+(i-1)/2$ to column $i$ if $i=1,3,5\ldots,n-1$; assign $3n+3+i/2$ if $i=2,4,6,\ldots,n-2$. We have used integers in $[3n+4,4n+2]$.
  \item In row 3, assign $3n+3-i$ to column $1\le i\le n-1$. We have used integers in $[2n+4,3n+2]$.
\end{enumerate} 

\nt  The resulting matrix is given by the following table:
\begin{center}
\nt Table 3: assignment of integers in $[2n+4,3n+2]\cup[3n+4,5n+1]$\\
$\setlength{\extrarowheight}{2pt}
\begin{array}{|c|c|c|c|c|c|c|c|c|}\hline
4n+3 & \frac{9n}{2}+3 & 4n+4 & \frac{9n}{2}+4 & \cdots & 5n & \frac{9n}{2}+1 & 5n+1 & \frac{9n}{2}+2 \\[2pt]\hline
\frac{7n}{2}+3 & 3n+4 & \frac{7n}{2}+4 & 3n+5 & \cdots &  \frac{7n}{2}+1 & 4n+1 &  \frac{7n}{2}+2  & 4n+2 \\[2pt]\hline
3n+2 & 3n+1 & 3n & 3n-1 & \cdots & 2n+7 & 2n+6 & 2n+5 & 2n+4\\\hline
\end{array}$
\end{center}

\ms\nt We next arrange integers in $[2,2n-1]\cup[2n+4,3n+2]$ as a $3\times n$ matrix as follows:
\begin{enumerate}[(1).]
  \item In row 1, assign $3n+3-i$ to column $1\le i\le n-1$. We have used integers in $[2n+4,3n+2]$.
  \item In row 2, assign $3n/2+(i-1)/2$ to column $i$ if $i=1,3,5,\ldots,n-1$; assign $n+i/2$ to column $i$ if $i=2,4,6,\ldots,n-2$. We have used integers in $[n+1,2n-1]$.
  \item In row 3, assign $(i+3)/2$ to column $i$ if $i=1,3,5,\ldots,n-1$; assign $n/2+1+i/2$ to column $i$ if $i=2,4,6,\ldots,n-2$. We have used integers in $[2,n]$.
\end{enumerate}

\nt  The resulting matrix is given by the following table:
\begin{center}
\nt Table 4: assignment of integers in $[2,2n-1]\cup[2n+4,3n+2]$\\
$\setlength{\extrarowheight}{2pt}
\begin{array}{|c|c|c|c|c|c|c|c|c|}\hline
3n+2 & 3n+1 & 3n & 3n-1 & \cdots & 2n+7 & 2n+6 & 2n+5 & 2n+4\\\hline
\frac{3n}{2} & n+1 & \frac{3n}{2}+1 & n+2 & \cdots & \frac{3n}{2}-2 & 2n & \frac{3n}{2}-1 & 2n-1 \\[2pt]\hline
2 & \frac{n}{2}+2 & 3 & \frac{n}{2}+3 & \cdots & n-1 &  \frac{n}{2}  & n &  \frac{n}{2}+1 \\[2pt]\hline
\end{array}$
\end{center}

\ms\nt For $1\le k\le 3$, $1\le i\le n-1$, let $c_{k,i}$ be the $(k,i)$-entry of Table 3, and $d_{k,i}$ be the $(k,i)$-entry of Table 4. Note that $d_{1,i} = c_{3,i}$. Define $f(u_iu^i_1) = c_{1,i}$, $f(xu_i) = c_{2,i}$, $f(u_iv_i) = c_{3,i}$, $f(xv_i) = d_{2,i}$ and $f(v_iv^i_1) = d_{3,i}$. It is obvious that $f$ is a bijective function.

\ms\nt Now, column sum of each column of Table 3 is $21n/2+8$. Thus, $w(u_i)=21n/2+8$ and $w(u^i_1)\in[4n+3,5n+1]$ for $1\le i\le n-1$. Similarly, the column sum of each column of Table 4 is $9n/2+4$. Thus, $w(v_i)=9n/2+4$ and $w(v^i_1)\in[2,n]$ for $1\le i\le n-1$. Moreover, $w(x) = [2n + (2n+1) + (3n+3)] + (3n+4) + \cdots + (4n+2) + (n+1) + \cdots + (2n-1) = (7n+4) + (n-1)(5n+3) = 5n^2+5n+1$.  Clearly,  for $1\le i\le n-1$, $w(x) \ne w(u^i_1) \ne  w(u_i) \ne w(v^i_1) \ne w(u^1_1) \ne w(v^1_1) \ne w(x_1)$. Note that $4n+3 \le w(u_n) = 4n+3 \ne w(v_n) = 4n+5 \le 5n+1$ for even $n\ge 6$. Therefore, $f$ is a local antimagic labeling that induces $2n+3$ distinct vertex colors. Consequently, $\chi_{la}(f_n \circ O_1) = 2n+3$ for even $n\ge 6$. \end{proof}

\begin{example} Figures~\ref{fig:n3m1} and~\ref{fig:n6m1} below give the labelings of $f_3\circ O_1$ and $f_6\circ O_1$ according to the proof in Theorem~\ref{thm-fnO1}.

\begin{figure}[h] 
  \centering
  \includegraphics[bb=6 12 369 223,width=3.5in,height=2in,keepaspectratio]{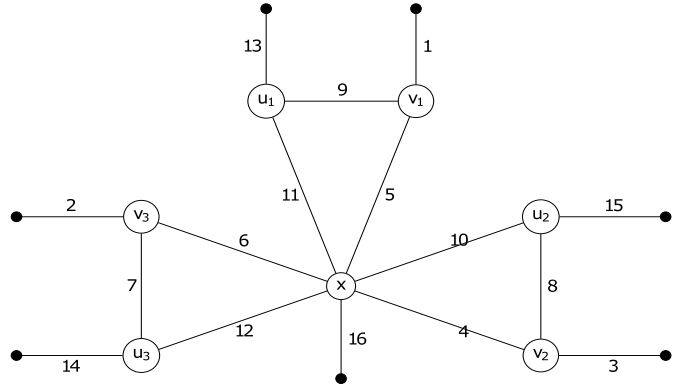}
  \caption{$\chi_{la}(f_3\circ O_1) = 9$ with induced vertex colors in $[1,3]\cup[13,16]\cup\{33,64\}$}
  \label{fig:n3m1}
\end{figure}

\begin{figure}[h] 
  \centering
  \includegraphics[bb=13 31 421 385,width=4.45in,height=3.5in,keepaspectratio]{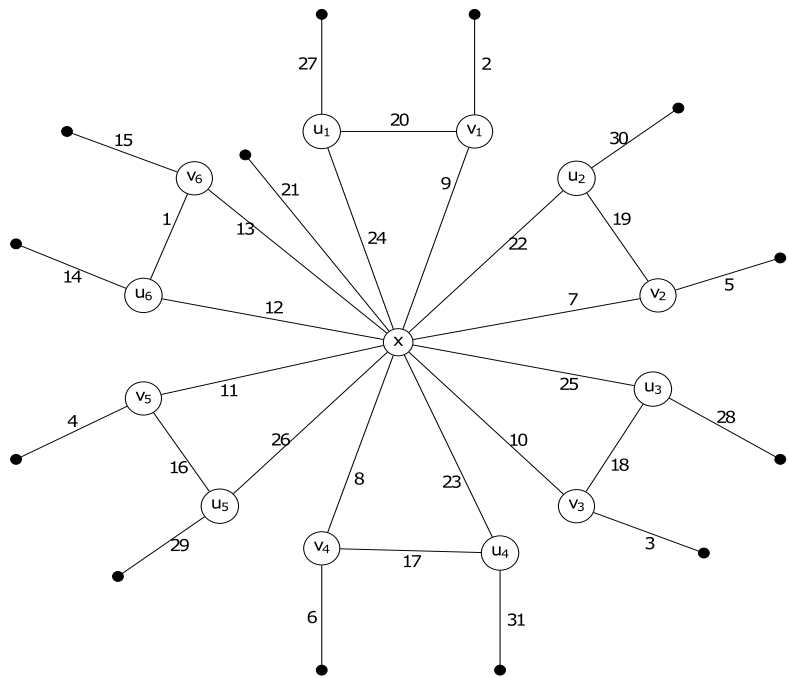}
  \caption{$\chi_{la}(f_6\circ O_1) = 15$ with induced vertex colors in $[1,6]\cup[27,31]\cup\{21,71,211\}$}
  \label{fig:n6m1}
\end{figure}


\end{example}


\end{document}